\documentclass{article}
\usepackage{latexsym}
\usepackage{amssymb}
\usepackage{amsfonts}
\usepackage{amsmath}
\usepackage{graphicx}
\usepackage{amsfonts}
\usepackage{epstopdf}
\usepackage{color}
\usepackage{enumerate}
\usepackage[latin1]{inputenc}
\usepackage{mathtools}
\usepackage{amsthm}
\usepackage[latin1]{inputenc}
\usepackage{mathrsfs}
\usepackage{graphicx}
\usepackage[latin1]{inputenc}
\newcommand{\R}{\mathbb{R}}

\newcommand{\fs}{{\cal F}_t}

\newcommand{\FF}{{\cal F}}
\newcommand{\De}{\Delta}
\newcommand{\prdi}{\prod_{i=1}^d (t_i - s_i )}

\def\lf{\left\lfloor}
\def\rf{\right\rfloor}

\newtheorem{teo}{Theorem}[section]

\begin{document}
\title{Strong approximations of 
Brownian sheet 
by uniform transport processes}
\date{}
\author{Xavier Bardina\footnote{X. Bardina is supported by the grantMTM2015-67802-P from SEIDI, Ministerio de Economia y Competividad.}, Marco Ferrante$^\ddagger$
 and Carles Rovira\footnote{C. Rovira is supported by the grant MTM2015-65092-P from MINECO/FEDER, UE.}}

\maketitle

$^*${\rm Departament de Matem\`atiques, Facultat de Ci\`encies,
Edifici C, Universitat Aut\`onoma de Barcelona, 08193 Bellaterra}.
{\tt Xavier.Bardina@uab.cat}
\newline
$\mbox{ }$\hspace{0.1cm} $^\ddagger${\rm Dipartimento di Matematica ``Tullio Levi-Civita",
Universit\`a degli Studi di Padova, Via Trieste 63, 35121 Padua, Italy}. {\tt marco.ferrante@unipd.it}
\newline
$\mbox{ }$\hspace{0.1cm} $^\dagger${\rm Departament de Matem\`atiques i Inform\`atica,
Universitat de Barcelona, Gran Via 585, 08007 Barcelona}. {\tt carles.rovira@ub.edu}

\begin{abstract}
Many years ago, Griego, Heath and Ruiz-Moncayo proved that it is possible to define
realizations of a sequence of uniform transform processes that converges almost surely
to the standard Brownian motion, uniformly on the unit time interval.
In this paper we extend their results to the multi parameter case. We begin
constructing a family of processes, starting from a set of independent standard Poisson processes,
that has realizations that converge
almost surely to the Brownian sheet, uniformly on the unit
square.
At the end the extension to the $d$-parameter Wiener processes is presented.

\end{abstract}

\section{Introduction}

Let $W=\{W(s,t): (s,t)\in [0,1]^{2}\}$ be a Brownian sheet, i.e. a zero mean real
continuous Gaussian process with covariance function
$E[W_{(s_1,t_1)} W_{(s_2,t_2)}]=(s_1\wedge s_2)(t_1\wedge t_2)$ for any
$(s_1,t_1), (s_2,t_2)\in [0,1]^2$.
The purpose of this paper is to find strong approximations of the Brownian sheet by
processes constructed from a family of independent standard Poisson processes.
In particular, we seek an extension, first in
the two-parameter case, of a result proved by Griego, Heath and Ruiz-Moncayo~\cite{art G-H-RM},
where the authors  present
realizations of a sequence of the uniform transform processes that converges almost surely
to the standard Brownian motion, uniformly on the unit time interval.
These results are not just interesting from a purely mathematical point of view, but are
of great interest in order to provide sound approximation strategies to solutions of stochastic
differential equations (in the classical Brownian motion case) and of
stochastic partial differential equations (in the present case of multi-parameter Wiener processes).

Griego, Heath and Ruiz-Moncayo  \cite{art G-H-RM}  deal with a sequence of uniform transport processes
\begin{equation}\label{aproxbm}
X_n(t)=\frac{1}{\sqrt{n}} (-1)^{A}\int_0^{tn}(-1)^{N(u)}du,\end{equation}
where $N=\{N(t), t\ge 0\}$ is a standard Poisson process and  $A$ is a random variable with law
\textrm{Bernoulli}$\left(\frac12\right)$  independent of  $N$.
To pass for the convergence in distribution, proved by Pinsky in \cite{pinsky1} to the almost sure convergence,
they make use of an embedding result due to Skorokhod (see \cite{sko}, page 163).

In the literature, there are some extensions of this result of almost sure convergence.
In \cite{art G-G2} Gorostiza and Griego
extended the result of \cite{BBR}) to the case of a diffusions. Again Gorostiza and Griego
\cite{art G-G} and Cs\"{o}rg\H{o} and Horv\'ath \cite{CH} obtained the
rate of convergence of the approximation sequence.
More recently, Garz\'on, Gorostiza and Le\'on \cite{GGL} defined a sequence of
processes that converges strongly to fractional Brownian motion
uniformly on bounded intervals, for any Hurst parameter $H\in(0,1)$
and computed the rate of convergence. In \cite{GGL2} and \cite{GGL3}
the same authors deal with subfractional Brownian motion and
fractional stochastic differential equations.
Finally in \cite{BBR}, Bardina, Binotto and Rovira proved the strong convergence to
a complex Brownian motion. As far as we know, our work is the first extension to the multiparameter case.

To the best of our knowledge, in the case of the Brownian sheet, or more generally of the
$d$--parameter Wiener processes, similar results have been proved just for the
weak convergence.
For instance, Bardina and Jolis \cite{BJ} prove that the process
\begin{equation*}
\frac{1}{n} \int_0^{tn}  \int_0^{sn} \sqrt{x y} (-1)^{N(x,y)}dxdy,
\end{equation*}
where $\{N(x,y), x\ge0, y \ge 0 \}$ is a Poisson process in the plane, converges in
law to a Brownian sheet when $n$ goes to infinity
and Bardina, Jolis and Rovira \cite{BJR} extended this result to get the weak
convergence to the $d$--parameter Wiener processes.

The present paper fills the gap in the case of a multi-parameter
Brownian motion. Indeed, in the next section we construct a uniform like transport process which
converges almost surely to the Brownian sheet and we present in the last section
the easy extension to the general case of a $d$--parameter Wiener process.

\section{Main result}

Let us built up our approximation processes.
Following some ideas of Bass and Pyke \cite{BP},
we will start defining a suitable partition of the unit square
For any $n$ and fixed $\lambda >0$, we can consider the partition of  the unit square
$[0,1]^2$ in disjoint rectangles
\[
\Big([0,\frac{1}{n^\lambda}] \times[0,1]\Big) \cup
\Big(  \bigcup_{k=2}^{\lf n^\lambda\rf}  (\frac{k-1}{n^\lambda} , \frac{k}{n^\lambda}] \times[0,1]  \Big)
\cup
\Big((\frac{\lf n^\lambda\rf}{n^\lambda} , 1] \times[0,1]\Big) .
\]
where $\lf x \rf$ denotes the greatest integer less than or equal to $x$.
If $W=\{W(s,t): $ $\,(s,t)\in [0,1]^{2}\}$ is a Brownian sheet on the unit square,
let $W^k$ denotes its restriction to each of the above defined rectangles
$ (\frac{k-1}{n^\lambda} , \frac{k}{n^\lambda}] \times [0,1]$.
That is, we define
\[
W^k(t):=W(\frac k{n^\lambda},t)-W(\frac{k-1}{n^\lambda},t),
\]
for $k\in\{1,2,\dots,\lf n^\lambda\,\rf\}$.
Thus, for any $l \in\{1,2,\dots,\lf n^\lambda\,\rf\}$ and $t \in [0,1]$
\[
W(\frac{l}{n^\lambda},t) = \sum_{k=1}^l W^k(t).
\]
Moreover, putting $\tilde W^k(t):=n^{\frac{\lambda}{2}} W^k(t)$, we obtain a
family
\[
\{\tilde W^k;\,k\in\{1,2,\ldots \lf n^\lambda\,\rf\} \}
\]
of
independent standard  Brownian motions defined in $[0,1]$.

From the paper of Griego, Heath and Ruiz-Moncayo  \cite{art G-H-RM} we know that
there exist realizations of  processes of  type (\ref{aproxbm})
that converge strongly and uniformly on bounded time intervals
to Brownian motion.
Using an approximation sequence for each one of the standard Brownian motions
$\{\tilde W^k;\,k\in\{1,2,\dots,\lf n^\lambda\,\rf\} \}$,  we will approximate the Brownian sheet by a process
$W_n$ such that for any   $l \in\{1,2,\dots,[n^\lambda\,]\}$ and $t \in [0,1]$
$$W_n(\frac{l}{n^\lambda},t) = \sum_{k=1}^l \frac{1}{n^\frac{1+\lambda} 2} (-1)^{A_k}\int_0^{nt}(-1)^{N_k(u)}du,
$$
where $\{N_k, k \ge 1\}$ is a family of independent standard Poisson processes and  $\{A_k, k \ge 1\}$ is a
sequence of independent random variables with law
$\textrm{Bernoulli}\left(\frac12\right)$, independent of the Poisson
processes.

Furthermore, using linear interpolation, we can define $W_n(s,t)$
on the whole unit square as follows:
\begin{eqnarray}
W_{n} (s,t)&=&\sum_{k=1}^{\lf s n^\lambda\rf}\frac1{n^{\frac{1+\lambda}{2}}}
(-1)^{A_k} \int_0^{tn}(-1)^{N_k(u)}du \label{pp2}
\\
&  & \qquad +(s n^\lambda -\lf s n^\lambda\rf)\frac1{n^{\frac{1+\lambda}{2}}}
(-1)^{A_{\lf s n^\lambda\rf+1}}\int_0^{tn}(-1)^{N_{\lf s n^\lambda\rf+1}(u)}du,
\nonumber
\end{eqnarray}
for any $(s,t) \in [0,\frac{\lf n^\lambda\rf}{n^\lambda} ]\times [0,1]$ and
$W_{n} (s,t)=W_{n} (\frac{\lf n^\lambda\rf}{n^\lambda},t)$ for any
$(s,t) \in [\frac{\lf n^\lambda\rf}{n^\lambda},1 ]\times [0,1]$.

\begin{center}
   \includegraphics[scale=0.5]{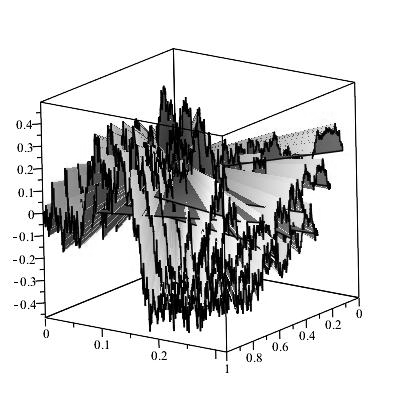}
   \end{center}
\begin{center}
   \begin{footnotesize}
  Simulation of $W_n$  with $\lambda:=\frac{1}{5.0001}$ and $n:=10^{\frac{1}{\lambda}}$.
   \end{footnotesize}
   \end{center}

The main result of this paper states as follows:

\begin{teo}\label{resultat}
There exists realizations of the process $\{ W_{n} (s,t), (s,t) \in [0,1]^2 \}$ with $\lambda \in (0,\frac15)$
on the same probability space as a  Brownian sheet $\{ W (s,t), (s,t) \in [0,1]^2 \}$
 such that
  $$ \lim_{n\rightarrow \infty}
 \max_{0\leq s,t\leq1}
 |W_n(s,t)-W(s,t)|=0 \quad a.s. $$
\end{teo}

\section{Proof of the main result}

The key of the proof of the convergence result in  \cite{art G-H-RM} is a Skorokhod's result (see \cite{sko} page 163) of a reproduction of independent random variables by evaluating Brownian motion at random times. In our proof, we will use the following obvious extension of Skorokhod's theorem:

\begin{teo}\label{skonostre}

Suppose that $\{\xi_1^k,\xi_2^k,\ldots,\xi_n^k; 1 \leq k \leq m \}$ are independent random variables such that
$E(\xi_i^k)=0$ and $Var(\xi_i^k)< \infty, 1 \leq i \leq n; 1 \leq k \leq m$ and that $\omega^1(t),\ldots,\omega^m(t)$
are a family of $m$ independent Brownian processes. Then, there exist nonegative independent random variables
$\tau_1^k,\tau_2^k,\ldots,\tau_n^k,  1 \leq k \leq m,$ for which the variables
\[
\omega^k(\tau_1^k), \omega^k(\tau_1^k+\tau_2^k)-\omega^k(\tau_1^k),\ldots,
\omega^k(\tau_1^k+\tau_2^k+\ldots+\tau_n^k)-
\omega(\tau_1^k+\tau_2^k+\ldots+\tau_{n-1}^k),
\]
have the same distribution as do $\xi_1^k, \xi_2^k,\ldots,\xi_n^k; 1 \leq k \leq m. $ Also:

(a) $E(\tau_i^k)=Var(\xi_i^k).$

(b) There exist a constant $L_p>0$ such that $E(\tau_i^k)^p \le L_p E(\xi_i^k)^{2p}.$

(c) For any $h>0$, if $\vert \xi_i^k \vert \le h,$ then
$
\vert \omega^k(s) - \omega^k(\sum_i^l \tau_i^k ) \vert \le h$ for $s \in [\sum_i^l\tau_i^k ,\sum_i^{l+1} \tau_i^k ].$

\end{teo}

\noindent
{\bf Proof of Theorem \ref{resultat}.}  The proof of our main result follows the method presented by Griego, Heath and Ruiz-Moncayo  \cite{art G-H-RM}.

For each $n$ and fixed $\lambda >0$, consider
$\{\eta_1^{(n)k},\eta_2^{(n)k},\ldots, \eta_{2n^2}^{(n)k}; 1 \le k \le \lf n^\lambda\rf \}$ families of
independent and identically distributed random variables with
exponential distribution of parameter $2n$, and
$\{K_1^{k},K_2^{k},\ldots, K_{2n^2}^{k}; 1 \le k \le \lf n^\lambda\rf \}$ sequences of independent and identically
distributed random variables so that
$P(K_i^{k}=1)=P(K_i^{k}=-1)=\frac12$.
Then, $\{K_1^{k}\eta_1^{(n)k},K_2^{k}\eta_2^{(n)k},\ldots, K_{2n^2}^{k}\eta_{2n^2}^{(n)k}; 1 \le k \le [n^\lambda] \}$ are
families of centered independent and identically distributed random
variables with common
$Var(K_i^{k}\eta_i^{(n)k})=\frac1{2n^2}$.

By Theorem \ref{skonostre}, the extension of the  Skorokhod embedding theorem,
for each $n\geq0$ and $1 \le k \le \lf n^\lambda\rf $ there
exist independent families of stopping times $\sigma_1^{(n)k},
\sigma_2^{(n)k},\ldots,\sigma_{2n^2}^{(n)k}$ that are nonnegative, independent and
identically distributed random variables such that the families
$$\tilde W^{k}(\sigma_1^{(n)k}),\tilde
W^{k}(\sigma_1^{(n)k}+\sigma_2^{(n)k}),\ldots, \tilde W^{k}(\sum_{i=1}^{2n^2} \sigma_i^{(n)k})$$ has the same law as
$$K_1^{k}\eta_1^{(n)k},K_1^{k}\eta_1^{(n)k}+K_2^{k}
\eta_2^{(n)k}, \ldots,  \sum_{i=1}^{2n^2}K_i^{k}\eta_i^{(n)k}.$$

Define now, for each $n\in {\mathbb N}$, $1 \le k \le \lf n^\lambda\rf $ and $i \in \{1,\ldots,2n^2\}$
$$
\gamma_i^{(n)k}=\frac1n |\tilde
W^{k}(\sum_{j=0}^i\sigma_j^{(n)k})-\tilde
W^{k}(\sum_{j=0}^{i-1}\sigma_j^{(n)k})|.
$$
Note that $\{ \gamma_i^{(n)k}: 1\le i \le 2n^2\}$ are families of independent and
identically distributed random variables with exponential
distribution of parameter $2n^2$.

We consider now the processes $\tilde W^{(n)k}(t),\,t\geq0,$ defined
by piecewise linear interpolation in such a manner that
\begin{equation}\tilde W^{(n)k}(\sum_{j=0}^i\gamma_j^{(n)k})=\tilde
W^{k}(\sum_{j=0}^i\sigma_j^{(n)k}), \label{bb1}\end{equation}
$\tilde W^{(n)k}(0)=0$.
Putting
$$ W^{(n)k}(t)=\frac1{n^\frac{\lambda}{2}}\tilde W^{(n)k}(t),$$
we can define the process
$$ W_n (\frac{l}{n^\lambda},t)=\sum_{k=1}^l W^{(n)k}(t)$$
for any $0\le l \le \lf n^\lambda\rf, t \in [0,1].$
Finally, we can define $W_n(s,t)$
for any $(s,t) \in [0,1]^2$ using again linear interpolation, that is,
if $s \in (\frac{k-1}{n^\lambda},\frac{k}{n^\lambda}]$,
$W_n(s,t)=(\frac{k}{n^\lambda}-s) W_n(\frac{k-1}{n^\lambda},t)+ ( s -  (\frac{k-1}{n^\lambda}))W_n(\frac{k}{n^\lambda},t)$
and $W_{n} (s,t)=W_{n} (\frac{\lf n^\lambda\rf}{n^\lambda},t)$ if $s \geq \frac{\lf n^\lambda\rf}{n^\lambda}$.
It is easy to check that  $W_n$ is a realization of the process (\ref{pp2})  .

Observe now that using Kolmogorov's inequality
\begin{eqnarray}
&&P(\max_{\footnotesize\begin{array}{cc}1\leq i\leq 2n^2\\1\leq
l\leq \lf n^\lambda\rf
\end{array}}|\gamma_1^{(n)l}+\cdots+\gamma_i^{(n)l}-\frac{i}{2n^2}|\geq
\varepsilon)\nonumber\\
&\leq& \sum_{1\leq l\leq \lf n^\lambda\rf}
P(\max_{\footnotesize\begin{array}{cc}1\leq i\leq 2n^2
\end{array}}|\gamma_1^{(n)l}+\cdots+\gamma_i^{(n)l}-\frac{i}{2n^2}|\geq
\varepsilon)\nonumber\\
&\leq& \frac{n^\lambda}{\varepsilon^2} \sum_{i=1}^{2n^2} Var( \gamma_i^{(n)l} ) =
\frac{n^\lambda}{2\varepsilon^2n^2}=\frac{1}{2\varepsilon^2n^{2-\lambda}}.\label{ppp}
\end{eqnarray}

By the same arguments and using that
$Var(\sigma_i^{(n)l})=\frac{3 L_2}{2 n^4}$, where $L_2$ is a constant obtained from Theorem \ref{skonostre}, we have also that
\begin{equation}
P(\max_{\footnotesize\begin{array}{cc}1\leq i\leq 2n^2\\1\leq
l\leq \lf n^\lambda \rf
\end{array}}|\sigma_1^{(n)l}+\cdots+\sigma_i^{(n)l}-\frac{i}{2n^2}|\geq
\varepsilon) \leq  \frac{3 L_2}{\varepsilon^2n^{2-\lambda}}.\label{ppp2}
\end{equation}

Since the Brownian sheet is $\alpha$--H\"older continuous for any $\alpha <\frac12$ and the
Lipschitz continuity of $W^{(n)}$,
 we obtain that
\begin{eqnarray}
&&\lim_{n\to\infty}\max_{\footnotesize\begin{array}{cc}0 \le s\le 1\\ 0 \le
t \le 1
\end{array}}|W^{(n)}(s,t)-W(s,t)|\nonumber\\
&=&\lim_{n\to\infty}\max_{\footnotesize\begin{array}{cc}1\leq i\leq
2n^2\\1\leq k\leq [n^\lambda\,]\end{array}} |W^{(n)}(\frac
k{n^\lambda},\frac
i{2n^2})-W(\frac k{n^\lambda},\frac i{2n^2})|\nonumber\\
&=&\lim_{n\to\infty}\max_{\footnotesize\begin{array}{cc}1\leq i\leq
2n^2\\1\leq k\leq [n^\lambda\,]\end{array}}
|\sum_{l=1}^kW^{(n)l}(\frac
i{2n^2})-\sum_{l=1}^kW^{l}(\frac
i{2n^2})|.\label{as}
\end{eqnarray}

Let us check that we can change $\sum_{l=1}^kW^l(\frac i{2n^2})$ by $\sum_{l=1}^kW^l(\sum_{j=0}^i\sigma_j^{(n)l})$ and
 $\sum_{l=1}^kW^{(n)l}(\frac i{2n^2})$ by $\sum_{l=1}^kW^{(n)l}(\sum_{j=0}^i\gamma_j^{(n)l})$.
Set
\[
A_{n,\delta}:=\{ \max_{\footnotesize\begin{array}{cc}1\leq i\leq 2n^2\\1\leq
l\leq [n^\lambda\,]
\end{array}}|\sigma_1^{(n)l}+\cdots+\sigma_i^{(n)l}-\frac{i}{2n^2}|\geq
\delta\};
\]
we get
\begin{eqnarray*}
&&P \Big( \max_{\footnotesize\begin{array}{cc}1\leq i\leq
2n^2\\1\leq k\leq \lf n^\lambda\rf\end{array}}
|\sum_{l=1}^kW^{l}(\frac
i{2n^2})-\sum_{l=1}^kW^l(\sum_{j=0}^i\sigma_j^{(n)l})| > \varepsilon \Big )
\\
&  \le & P(A_{n,\delta})  +
P \Big( \{ \max_{\footnotesize\begin{array}{cc}1\leq i\leq
2n^2\\1\leq k\leq \lf n^\lambda\rf \end{array}}
|\sum_{l=1}^kW^{l}(\frac
i{2n^2})-\sum_{l=1}^kW^l(\sum_{j=0}^i\sigma_j^{(n)l})| > \varepsilon \}  \cap
A_{n,\delta}^c \Big) \\ &  \le &
\frac{3 L_2}{\delta^2n^{2-\lambda}}
+
P \Big( \max_{\footnotesize\begin{array}{c}1\leq k\leq [n^\lambda\,]\end{array}}
\sum_{l=1}^k   \sup_{\footnotesize\begin{array}{c}\vert s-t\vert < \delta\end{array}}
\vert W^{l}(t) - W^l(s) \vert  > \varepsilon  \Big) \\ &  \le &
\frac{3 L_2}{\delta^2n^{2-\lambda}}
+  n^\lambda
P \Big(
  \sup_{\footnotesize\begin{array}{c}\vert s-t\vert < \delta \end{array}}
   \vert W^{1}(t) - W^1(s) \vert  > \frac{\varepsilon}{n^\lambda}  \Big)
\\ &  \le &
\frac{3 L_2}{\delta^2n^{2-\lambda}}
+ \frac{ n^{\lambda(M+1)}}{\varepsilon^M}
E \Big(
  \sup_{\footnotesize\begin{array}{c}\vert s-t\vert < \delta \end{array}}
  \vert W^{1}(t) - W^1(s) \vert^M  \Big)
\\ &  \le &
\frac{3 L_2}{\delta^2n^{2-\lambda}}
+ C_M \frac{ n^{\lambda(M+1)}}{\varepsilon^M}
 \Big(  \delta \log \big( \frac{2}{\delta} \big) \Big)^\frac{M}{2},
 \end{eqnarray*}
where we have used (\ref{ppp2}), Markov's inequality and estimates for the modulus of continuity of the Brownian motion.
Putting $\delta=n^{\frac{\lambda-1}{2} + \mu}$, for any $\lambda \in (0,\frac15)$, we can choose $\mu>0$  such that
 $5 \lambda - 1 + 2\mu <0$. Then
\[
\frac{\lambda-1}{2} + \mu <0
\]
and there exists $M$
large enough such that
\[
\lambda(M+1)+ \frac{M}{2} \big(\frac{\lambda-1}{2} + \mu \big)  = \frac{M}{4} \big( 5\lambda - 1 + 2\mu \big) + \lambda< -1.
\]
Thus, we  get that
\[
\sum_{n=1}^\infty
P \Big( \max_{\footnotesize\begin{array}{cc}1\leq i\leq
2n^2\\1\leq k\leq \lf n^\lambda\rf\end{array}}
|\sum_{l=1}^kW^{l}(\frac
i{2n^2})-\sum_{l=1}^kW^l(\sum_{j=0}^i\sigma_j^{(n)l})| > \varepsilon \Big ) < \infty.
\]
Then by the Borel-Cantelli lemma we have
\begin{equation}
\lim_{n\to\infty}\max_{\footnotesize\begin{array}{cc}1\leq i\leq
2n^2\\1\leq k\leq \lf n^\lambda\rf\end{array}}
|\sum_{l=1}^kW^{l}(\frac
i{2n^2})-\sum_{l=1}^kW^l(\sum_{j=0}^i\sigma_j^{(n)l})|=0, \qquad a.s. \label{as1}
\end{equation}
Using (\ref{ppp}) and similar argument we get
\begin{equation}
\lim_{n\to\infty}\max_{\footnotesize\begin{array}{cc}1\leq i\leq
2n^2\\1\leq k\leq \lf n^\lambda\rf\end{array}}
|\sum_{l=1}^kW^{(n)l}(\frac
i{2n^2})-\sum_{l=1}^kW^{(n)l}(\sum_{j=0}^i\gamma_j^{(n)l})|=0, \qquad a.s. \label{as2}
\end{equation}

\noindent
Putting together (\ref{as}),  (\ref{as1}) and  (\ref{as2}) and using  that  $\tilde W^k(t):=n^{\frac{\lambda}{2}} W^k(t)$  and
(\ref{bb1}), at the end we obtain
\begin{eqnarray*}
&&\lim_{n\to\infty}\max_{\footnotesize\begin{array}{cc}0 \le s\le 1\\ 0 \le
t \le 1
\end{array}}|W^{(n)}(s,t)-W(s,t)|\\&=&\lim_{n\to\infty}\max_{\footnotesize\begin{array}{cc}1\leq i\leq
2n^2\\1\leq k\leq \lf n^\lambda\rf\end{array}}
|\sum_{l=1}^kW^{(n)l}(\sum_{j=0}^i\gamma_j^{(n)l})-\sum_{l=1}^kW^l(\sum_{j=0}^i\sigma_j^{(n)l})|
\\
&=&\lim_{n\to\infty}\max_{\footnotesize\begin{array}{cc}1\leq i\leq
2n^2\\1\leq k\leq \lf n^\lambda\rf\end{array}}
|\frac1{n^{\frac{\lambda}{2}}}\sum_{l=1}^k\tilde W^{l}(\sum_{j=0}^i\sigma_j^{(n)l})-\sum_{l=1}^kW^l(\sum_{j=0}^i\sigma_j^{(n)l})|
\\
&=&0,
\end{eqnarray*}
a.s. and the proof is complete.
\hfill$\square$

\section{Extension to the $d$-parameter Wiener process}

Suppose $d\ge 2$ and consider $[0,1]^d  \subset \R^d$.
Define on $[0,1]^d$ the partial ordering
\[
s^1\le s^2 \Leftrightarrow
s^1_i \le s^2_i \ \mbox{for any} \ 1\le i \le d
\]
and denote $I_0=\{s\in [0,1]^d: s_1\cdot\cdots\cdot
s_d=0\}$.
Let $(\Omega,\FF,Q)$ be a complete probability space and let $\{\FF_s, s \in [0,1]^d\}$
be a family of sub-$\sigma$-fields of $\FF$ such that: $\FF_{s^1} \subseteq \FF_{s^2}$
for any $s^1 \le s^2$.
Fix $t \in [0,1]^d$ we also consider
$\FF^T_t := \bigvee_{i=1}^d \FF_{(1,\cdots,1,t_i,1,\cdots,1)}$ (the $\sigma$-field generated by all the past of $t$).
Given $s<t$ we denote by $\De_s X(t)$ the increment of the process
$X$ over the rectangle
$(s,t]=\prod_{i=1}^d (s_i, t_i] \subset \R^d.$

Then, a $d$-parameter continuous process $W= \{W(s); s \in [0,1]^d \subset \R^d_+ \}$
is called a $d$-parameter $\{ \fs \}$-Wiener process if it is $\{ \fs \}$-adapted,
null on $I_0$ and for any $s < t$ the increment $\De_s W(t)$ is
 independent of $\mathcal F_s^T$ and is normally distributed with zero
mean and variance $\prdi.$

In order to built up the approximation sequence we can deal now with the disjoint rectangles
\[
\prod_{i=1}^{d-1}   (\frac{k_i-1}{n^\lambda} , \frac{k_i}{n^\lambda}] \times[0,1]
\]
where $1 \le k_1,\ldots,k_{d-1} \le {\lf n^\lambda \rf}$
and the convention that the intervals are left-closed at the zero end points.
Set $W^{k_1,\ldots,k_{d-1}}$ the restriction of the $d$-parameter Wiener process  to each of these rectangles,
that is, we define
\[
W^{k_1,\ldots,k_{d-1}}(t):=\De_{(\frac{k_1-1}{n^\lambda},\ldots,\frac{k_{d-1}-1}{n^\lambda}, t)}
W(\frac{k_1}{n^\lambda},\ldots,\frac{k_{d-1}}{n^\lambda}, t),
\]
for  $1 \le k_1,\ldots,k_{d-1} \le {\lf n^\lambda\rf}$.
Thus, for any $1 \le l_1,\ldots,l_{d-1} \le {\lf n^\lambda \rf}$ and $t \in [0,1]$
\[
W(\frac{l_1}{n^\lambda},\ldots,\frac{l_{d-1}}{n^\lambda},t) =
\sum_{k_1=1}^{l_1}\ldots  \sum_{k_{d-1}=1}^{l_{d-1}} W^{k_1,\ldots,k_{d-1}}(t)
\]
and thus
$\tilde W^{k_1,\ldots,k_{d-1}}(t):=n^{(d-1)\frac{\lambda}{2}} W^{k_1,\ldots,k_{d-1}}(t)$ is a family of
independent standard  Brownian motions defined in $[0,1]$.
As in the 2-dimensional case we will approximate the $d$-dimensional Wiener process by a  process
$W_n$ such that for any   $1 \le l_1,\ldots,l_{d-1} \le {\lf n^\lambda\rf}$ and $t \in [0,1]$
\begin{eqnarray*}
&&W_n(\frac{l_1}{n^\lambda},\ldots,\frac{l_{d-1}}{n^\lambda},t)  \\
&&= \sum_{k_1=1}^{l_1}\ldots  \sum_{k_{d-1}=1}^{l_{d-1}}
\frac{1}{n^\frac{1+(d-1)\lambda} 2} (-1)^{A_{k_1,\ldots,k_{d-1}}}\int_0^{nt}(-1)^{N_{k_1,\ldots,k_{d-1}}(u)}du,
\end{eqnarray*}
where $\{N_{k_1,\ldots,k_{d-1}} \}$ is a family of independent standard Poisson processes and
$\{A_{k_1,\ldots,k_{d-1}} \}$ is a sequence of independent random variables with law
$\textrm{Bernoulli}\left(\frac12\right)$  independent of the Poisson
processes.
We can define $W_n(s_1,\ldots,s_{d-1},t)$ for any $(s_1,\ldots,s_{d-1}) \in [0,1]^{d-1}$
using a $d-1$-dimensional linear interpolation.

The main theorem reads as follows:

\begin{teo}\label{resultatd}
There exists realizations of the process $\{ W_{n} (s_1,\ldots,s_d):$ $(s_1,\ldots,$ $s_d) \in [0,1]^d \}$
with $\lambda \in (0,\frac{1}{5(d-1)})$
on the same probability space as a $d$-dimensional Wiener process
$\{ W(s_1,\ldots,s_d): (s_1,\ldots,s_d) \in [0,1]^d \}$
such that
\[
\lim_{n\rightarrow \infty} \max_{0\leq s_1,\ldots,s_d \leq1 } |W_n(s_1,\ldots,s_d)-W(s_1,\ldots,s_d)|=0 \quad a.s.
\]
\end{teo}

\noindent
{\it Sketch of the Proof:} the proof follows the same steps that for the Brownian sheet.
For each $n\geq0$ and $1 \le{ k_1,\ldots,k_{d-1}} \le \lf n^\lambda\rf $ there
exist independent families of stopping times $\sigma_1^{(n){ k_1,\ldots,k_{d-1}}},
\ldots,\sigma_{2n^2}^{(n){ k_1,\ldots,k_{d-1}}}$ that are nonnegative, independent and
identically distributed random variables that we can use to define, for each $n$ and
$1 \le { k_1,\ldots,k_{d-1}}\le \lf n^\lambda\rf $,
\[
\gamma_i^{(n){ k_1,\ldots,k_{d-1}}}=\frac1n |\tilde
W^{{ k_1,\ldots,k_{d-1}}}(\sum_{j=0}^i\sigma_j^{(n){ k_1,\ldots,k_{d-1}}})-\tilde
W^{{ k_1,\ldots,k_{d-1}}}(\sum_{j=0}^{i-1}\sigma_j^{(n){ k_1,\ldots,k_{d-1}}})|,
\]
for any $i \in \{1,\ldots,2n^2\}$.

We consider now the processes $\tilde W^{(n){ k_1,\ldots,k_{d-1}}}(t),\,t\geq0,$ defined
by piecewise linear interpolation in such a manner that
\begin{equation*}\tilde W^{(n){ k_1,\ldots,k_{d-1}}}(\sum_{j=0}^i\gamma_j^{(n){ k_1,\ldots,k_{d-1}}})=\tilde
W^{{ k_1,\ldots,k_{d-1}}}(\sum_{j=0}^i\sigma_j^{(n){ k_1,\ldots,k_{d-1}}}) \end{equation*}
$\tilde W^{(n){ k_1,\ldots,k_{d-1}}}(0)=0$.
Then we can define the process
\[
W^{(n)}(\frac{l_1}{n^\lambda},\ldots,\frac{l_{d-1}}{n^\lambda},t)=
\sum_{k_1=1}^{l_1}\ldots  \sum_{k_{d-1}=1}^{l_{d-1}} W^{k_1,\ldots,k_{d-1}}(t)
\]
for any $l_1,\ldots,l_{d-1} \in [0,\lf n^\lambda\rf]^{d-1}, t \in [0,1].$ Finally, we can define $W^{(n)}(s,t)$
for any $(s,t) \in [0,1]^2$ using $d-1$-dimensional linear interpolation.

Following the Brownian sheet case we get the following estimates:
\begin{eqnarray*}
& &P(\max_{\footnotesize\begin{array}{cc}1\leq i\leq 2n^2\\1\leq
l_1,\ldots,l_{d-1}\leq \lf n^\lambda\rf
\end{array}}|\gamma_1^{(n)l_1,\ldots,l_{d-1}}+\cdots+\gamma_i^{(n)l_1,\ldots,l_{d-1}}-\frac{i}{2n^2}|\geq
\varepsilon) \\ & & \qquad\qquad  \le \frac{1}{2\varepsilon^2n^{2-\lambda(d-1)}}\end{eqnarray*}
and
\begin{eqnarray*}
& &P(\max_{\footnotesize\begin{array}{cc}1\leq i\leq 2n^2\\1\leq
l_1,\ldots,l_{d-1}\leq \lf n^\lambda\rf
\end{array}}|\sigma_1^{(n)l_1,\ldots,l_{d-1}}+\cdots+\sigma_i^{(n)l_1,\ldots,l_{d-1}}-\frac{i}{2n^2}|\geq
\varepsilon) \\ & & \qquad\qquad \le \frac{3
L_2}{\varepsilon^2n^{2-\lambda(d-1)}}\end{eqnarray*} and
\begin{eqnarray*}
&&P \Big( \max_{\footnotesize\begin{array}{cc}1\leq i \leq
2n^2\\1\leq k_1,\ldots,k_{d-1}\leq \lf n^\lambda\rf\end{array}}
|\sum_{l_1=1}^{k_1}\ldots  \sum_{l_{d-1}=1}^{k_{d-1}} W^{l_1,\ldots,l_{d-1}}(\frac
i{2n^2} )
\\
& & \qquad -\sum_{l_1=1}^{k_1}\ldots  \sum_{l_{d-1}=1}^{k_{d-1}} W^{l_1,\ldots,l_{d-1}}(\sum_{j=0}^i\sigma_j^{(n)l})| > \varepsilon \Big )
\\
& \le &
\frac{3 L_2}{\delta^2n^{2-\lambda(d-1)}} \\ & & \quad
+
P \Big( \max_{\footnotesize\begin{array}{c}1\leq k_1,\ldots,k_{d-1}\leq \lf n^\lambda\rf\end{array}}
\sum_{l_1=1}^{k_1}\ldots  \sum_{l_{d-1}=1}^{k_{d-1}}  \\ & & \qquad  \sup_{\footnotesize\begin{array}{c}\vert s-t\vert < \delta\end{array}}    \vert W^{l_1,\ldots,l_{d-1}}(t) - W^{l_1,\ldots,l_{d-1}}(s) \vert  > \varepsilon  \Big)
\\
& \le & \frac{3 L_2}{\delta^2n^{2-\lambda}} +  n^{\lambda(d-1)} P
\Big( \sup_{\footnotesize\begin{array}{c}\vert s-t\vert < \delta
\end{array}} \vert W^{l}(t) - W^l(s) \vert  >
\frac{\varepsilon}{n^{\lambda(d-1)}}  \Big)
\\
& \le &
\frac{3 L_2}{\delta^2n^{2-\lambda(d-1)}}
+ C_M \frac{ n^{\lambda(d-1)(M+1)}}{\varepsilon^M}
 \Big(  \delta \log \big( \frac{2}{\delta} \big) \Big)^\frac{M}{2}.
 \end{eqnarray*}
We finish the proof, as in the Brownian sheet case,
putting $\delta=n^{\frac{\lambda(d-1)-1}{2} + \mu}$, for any $\lambda  \in (0,\frac{1}{5(d-1)})$, we can choose $\mu>0$  such that
 $5 \lambda (d-1)- 1 + 2\mu <0$. Then
\[
\frac{\lambda(d-1)-1}{2} + \mu <0
\]
and there exists $M$
large enough such that
\[
\lambda(M+1)+ \frac{M}{2} \big(\frac{\lambda(d-1)-1}{2} + \mu \big)  = \frac{M}{4} \big( 5\lambda(d-1) - 1 + 2\mu \big) + \lambda(d-1)< -1.
\]
\hfill$\square$

\section{Conclusions and future work}

In the present paper we have constructed a uniform transport like process that
converges strongly to the multi-parameter Wiener process.
Our result makes use of the 1-parameter result proved by
Griego, Heath and Ruiz-Moncayo~\cite{art G-H-RM}
and some ideas from
of Bass and Pyke \cite{BP}.

The interest of this class of results is not only purely mathematical, but also
for their possible application in order to approximate the solutions to stochastic
partial differential equations, which arise in many fields as physics, biology or finance.
We plan to continue our investigation with the aim to extend the present result to
these class of processes and to obtain results on the rate of convergence, that in
the present case looks quite challenging, needed for the applications.

\end{document}